\documentclass[authoryear, preprint,12pt]{elsarticle}

\usepackage{amssymb}
\usepackage{url}
\usepackage{graphicx}
\usepackage{amsmath}
\usepackage[gen]{eurosym}
\usepackage{wrapfig}
\usepackage{tabularx,ragged2e}
\usepackage{adjustbox}
\usepackage{booktabs}
\usepackage{tikz}
\usepackage{subfig}

\usepackage{enumitem}
\usepackage[gen]{eurosym}
\usepackage{algorithm}
\usepackage{algpseudocode}
\usepackage{enumerate}
\usepackage{threeparttable}
\usepackage{footnote}
\usepackage{multirow}

\usepackage{pgfplots}
\pgfplotsset{compat=1.11}
\usepgfplotslibrary{fillbetween}
\usetikzlibrary{intersections}

\newcommand{\rhs}{e}

\newcommand{\changes}[1]{\textcolor{black}{#1}}

\newtheorem{thm}{Theorem}

\newtheorem{rem}[thm]{Remark}

\newproof{pf}{Proof}

\usepackage{lipsum}
\makeatletter
\def\ps@pprintTitle{%
 \let\@oddhead\@empty
 \let\@evenhead\@empty
 \def\@oddfoot{}%
 \let\@evenfoot\@oddfoot}
\makeatother
\begin{document}

\begin{frontmatter}



\title{Multicriteria Adjustable Regret Robust Optimization for Building Energy Supply Design}


\author[label1]{Elisabeth Halser\corref{cor1}}
\ead{elisabeth.halser@itwm.fraunhofer.de}
\author[label1]{Elisabeth Finhold}
\ead{elisabeth.finhold@itwm.fraunhofer.de}
\author[label1]{Neele Leith\"auser}
\ead{neele.leithaeuser@itwm.fraunhofer.de}
\author[label1]{Tobias Seidel}
\ead{tobias.seidel@itwm.fraunhofer.de}
\author[label1]{Karl\nobreakdash-Heinz K\"ufer}
\ead{karl-heinz.kuefer@itwm.fraunhofer.de}

\cortext[cor1]{Corresponding Author}

\affiliation[label1]{organization={Fraunhofer Institute for Industrial Mathematics ITWM},
	addressline={Fraunhofer-Platz~1}, 
	city={Kaiserslautern},
	postcode={67663}, 
	country={Germany}}

\begin{abstract}
Optimizing a building's energy supply design is a task with multiple competing criteria, where not only monetary but also, for example, an environmental objective shall be taken into account. Moreover, when deciding which storages and heating and cooling units to purchase (here-and-now-decisions), there is uncertainty about future developments of prices for energy, e.g. electricity and gas. This can be accounted for later by operating the units accordingly (wait-and-see-decisions), once the uncertainty revealed itself. Therefore, the problem can be modeled as an adjustable robust optimization problem. We combine adjustable robustness and multicriteria optimization for the case of building energy supply design and solve the resulting problem using a column and constraint generation algorithm in combination with an $\varepsilon$-constraint approach.

In the multicriteria adjustable robust problem, we simultaneously minimize worst-case cost regret and carbon emissions. We take into account future price uncertainties and consider the results in the light of information gap decision theory to find a trade-off between security against price fluctuations and over-conservatism. We present the model, a solution strategy and discuss different application scenarios for a case study building.
\end{abstract}

\begin{graphicalabstract}
\includegraphics{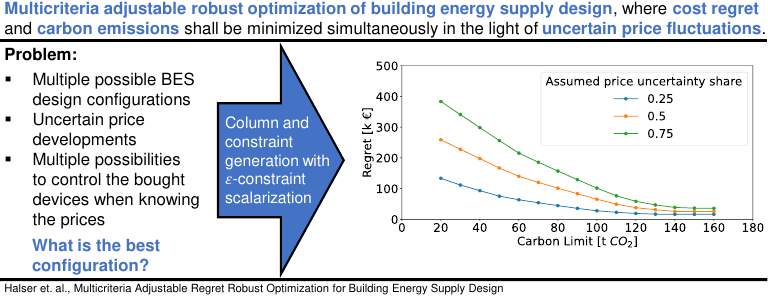}
\end{graphicalabstract}

\begin{highlights}
\item Novel multicriteria optimization model for building energy supply design
\item Consideration of cost regret with respect to price uncertainty as objective
\item Column and constraint generation algorithm and proof of convergence
\item Demonstration of the approach in a case study for an office building at an early planning stage
\end{highlights}


\end{frontmatter}



\section{Introduction}
Buildings and their energy consumption are responsible for about 40\% of global greenhouse gas emissions \citep{IEA} and for a significant amount of many companies' annual costs. In the following, we use the term carbon emissions to consider all types of greenhouse gas emissions as they can be converted to carbon equivalent emissions. The two possibly competing criteria of minimizing costs and carbon emissions make the task of designing a building's energy supply an inherently multicriteria one for which trade-offs have to be identified. 

The existence of multiple criteria is not the only challenge in building energy supply design. A second structural difficulty is that the problem naturally comes with different decision stages. Initially, there is the first stage device purchasing decision (\emph{here-and-now} decision). Throughout the article, we use the term device to talk about all types of heating and cooling units and storages as well as photovoltaic. 
Before the actual operation of the purchased units, the costs for operating them, which are highly uncertain due to potential price fluctuations, become known. 
Then, in a second decision stage, there is the possibility to react to the then observed resource price reality (\emph{wait-and-see} decisions) by controlling the heating and cooling units. The reason for distinguishing between here-and-now and wait-and-see variables is that at the beginning, the prices are uncertain. Because of the two stages, where the second stage can be adjusted to the realized uncertainty, we call the problem \emph{adjustable}.

One way to deal with uncertainty is to consider the worst case of an assumed set of possible scenarios, which we call uncertainty set. This approach is called robust optimization. When considering the worst possible price scenario, we call the problem structure multicriteria adjustable robust.

However, when solely focusing on the optimization of worst-case total costs, it is likely that the purchasing decision is tailored for an unlikely extreme scenario \changes{and} is likely to underperform in most other, more moderate, scenarios. Therefore, it is preferable to take a here-and-now decision which will cause minimal regret when the actual resource prices are known. Regret can be quantified as the difference between the minimal costs given a fixed purchasing decision and the minimal costs that would have been possible, if the prices had been known in advance. The concept for this kind of problem is called regret robust optimization. Hence, in this study, we employ cost regret alongside carbon emissions as an objective in our multicriteria adjustable regret robust optimization. 

Modeling building energy supply design as optimization problem is a large field of research. Most of the models are linear programs (LP) or mixed integer programs (MILP). Detailed modeling considerations are given by \cite{wirtz2021design}. 

A good overview over the mathematical field of multicriteria optimization (MCO) is for example given in the book by \cite{ehrgott2005multicriteria}. MCO problems have more than one objective. If some of these objectives are contradicting, it is not possible to minimize them at the same time. Therefore, compromises, also called trade-offs, have to be found. The basic idea of MCO is to identify solutions where no other solution outperforms them in every objective. Such solutions are called efficient, and their image in objective space is known as the Pareto front. One approach to approximate the Pareto front is the $\varepsilon$-constraint method, which inspired our approach. Multicriteria optimization has been applied to building energy supply design in different articles, for example in the work of \cite{richarz2020multicriteria} or \cite{finke2023implementing}.

Robust optimization is an approach to deal with uncertainty in optimization problems. A good overview is given by \cite{ben2009robust} or more recently by \cite{bertsimas2022robust}, where both books also cover adjustable robustness. Robust optimization for energy supply design was for example done by \cite{moret2020decision} and for building renovation by \cite{richarz2021robust}.

Regret robustness aims for a solution that does not minimize the worst-case outcome but the worst regret over all scenarios. Regret robustness ensures that the optimal solution is not only determined by a potentially unlikely extreme scenario, but rather by how much one would regret a decision when comparing it to what would have been possible if the scenario had been known in advance.
It was applied to the energy problem by Yokoyama et al. for continuous \citep{yokoyama2014revised} and discrete \citep{yokoyama2021hierarchical} design variables, where in these works the uncertainty comes from uncertain loads.

Robust multicriteria optimization allows for different definitions of efficient solutions, of which some are given by \cite{botte2019dominance}. The most common definition is that of point-based min-max robust efficiency, which for building energy supply design was applied by \cite{majewski2017robust}.

Multicriteria adjustable robust optimization (MARO) is a relatively new topic, that combines adjustable and multicriteria optimization and was studied in the linear decision rule case by \cite{chuong2022adjustable} and for general cases in \citep{halsermaro}. There, different computational approaches are compared, especially a weighted-sum scalarization and an $\varepsilon$-constraint scalarization. Solutions which are found with these approaches are called weighted sum MARO-efficient and constraint MARO-efficient, respectively. \cite{yokoyama2001multiobjective} used a weighted sum approach to deal with regret robustness against uncertain loads, which correspond to right-hand side uncertainty.

Inverse robustness~\citep{berthold2024unified}, also called information gap decision theory, was already applied in the energy context for one \citep{javadi2017robust} and multiple \citep{feng2023multi} objectives. It is a way of answering the question of what size of uncertainty set the solution should be robust against, which forms an additional objective on its own.

In this paper, we compute constraint MARO-efficient solutions for the problem of multicriteria building energy supply design with the objectives carbon emission and cost regret minimization. Constraint MARO-efficient solutions are beneficial for the problem as they come with good bounds for the objectives, which hold over all scenarios of the uncertainty set, and their solution in the objective space can be interpreted similarly to a usual Pareto front \citep{halsermaro}. Both properties are especially valuable for decision makers and make this approach valuable in practice. We will focus on uncertainty sets that consist of possible future price fluctuations. Compared to load uncertainties, price fluctuations are much more unpredictable and therefore hold practical relevance. In this paper, we present a case study in which we place the results in the context of inverse robustness by varying the size of the uncertainty set. In this way, the decision makers can weigh up how much price uncertainty they want to hedge against and what this might cost.

\begin{table}
\centering
\adjustbox{angle=90, width=0.42\textwidth}{
\begin{threeparttable}
		\begin{tabular}{lp{2.2cm}p{2.2cm}p{2.2cm}p{2.2cm}p{2.8cm}p{2.2cm}}
		& \cite{yokoyama2014revised} & \cite{feng2023multi} & \cite{yang2024two} & \cite{fan2014min} & \cite{aghamohamadi2019bidding} & our work \\
		\hline
		multicriteria optimization & X & $\surd$\tnote{ws} & $\surd$\tnote{ws}& X & X & $\surd$  \\
		adjustable & $\surd$ & X & $\surd$ & X & $\surd$ & $\surd$ \\
		price uncertainty & X & X & X & $\surd $ & $\surd$ & $\surd$ \\
		regret consideration & $\surd$ & X & X & $\surd$ & X & $\surd$ \\
		inverse robust idea & $\surd$\tnote{e} & $\surd$\tnote{o} & $\surd$\tnote{e} & X & $\surd$\tnote{e} & $\surd$\tnote{e} \\
		constraint generation & $\surd$ & X & $\surd$\tnote{ccg} & $\surd$ & $\surd$\tnote{ccg} & $\surd$\tnote{ccg}\\
		energy system design & $\surd$ & X & $\surd$ & X & X & $\surd$ \\ 
		\hline
	\end{tabular}
\begin{tablenotes}\footnotesize
	\item[ws] weighted sum scalarization without aim to explore different trade-offs in objective space
	\item[e] varying the sizes of the uncertainty set
	\item[o] maximizing the size of the uncertainty set
	\item[ccg] column and constraint generation
\end{tablenotes}
\end{threeparttable}}
	\caption{Classification of the contribution of this work.} \label{classification}
\end{table}

Table~\ref{classification} puts our work in the context of the already existing literature. There is a great number of further articles which could be added to the table as they deal with some of the topics. However, to the best of our knowledge, there has not been a combination of all these points so far.
To summarize our contributions, 
\begin{itemize}
	\item we compute constraint MARO-efficient solutions for the problem of minimizing cost regret and carbon emissions of building energy supply design, thereby combining regret robust optimization and multicriteria optimization.
	\item we consider uncertain prices in the objective function (instead of uncertain loads on the right-hand side as in \citep{yokoyama2014revised}) in the cost regret optimization and show how the algorithmic ideas of \cite{zeng2013solving} carry over.
\end{itemize}

The remainder of this article is structured as follows. In Section~\ref{model}, we give a brief overview of our model, for which the details are given in \ref{problemformulation}. In Section~\ref{maro} we explain how we apply the ideas of MARO to our problem. The detailed formulation of the solution algorithm is given in \ref{solutionalgorithm}. In Section~\ref{casestudy} we discuss a case study and summarize our results in Section~\ref{conclusion}.

\section{Model} \label{model}

As a basic model, yet without consideration of regret, we use a slightly adapted model from \citep{halserpn}, which we summarize here, while the detailed formulas are given in \ref{problemformulation}. Unlike as in \citep{halserpn}, we do not consider a third objective function of inconvenience but assume that all predetermined loads have to be fulfilled. 

The main goal in building energy supply design is to find a design configuration $\mathbf{d} \in \mathbb{R}^{n_d}, n_d \in \mathbb{N}$. For a list of commercially available unit types and storages, the vector $\mathbf{d}$ contains the selected size. In order to find $\mathbf{d}$ and evaluate annual costs and emissions, it is necessary to also find controls $\mathbf{s} \in \mathbb{R}^{n_s}, n_s \in \mathbb{N}$ for a representative period of time. A suitable set of representative days $\mathcal{U}$, which can be used to approximate annual costs and emissions, is found by k-Medoids clustering as suggested in \citep{schutz2018comparison}. For all representative days with all their time steps, the vector $\mathbf{s}$ represents the control states for all unit types and the storage states. 
The decision variables $\mathbf{d}$ and $\mathbf{s}$ have to fulfill the constraints that for all representative days the loads are fulfilled and the technical constraints are met. 
The objective expressions $costs^{p}(\mathbf{d},\mathbf{s})$ and $co_2(\mathbf{d},\mathbf{s})$, which are to be minimized, represent the annual costs and carbon emissions, respectively, caused by the design configuration $\mathbf{d}$ with controls $\mathbf{s}$ under price vector $p$. As we will need this split later, we introduce
\begin{align*}
	co_2(\mathbf{d}, \mathbf{s})= co_2^{inv}(\mathbf{d}) + co_2^{op}(\mathbf{d}, \mathbf{s}),
\end{align*}
where $co_2^{inv}(\mathbf{d})$ denotes the investment or life cycle carbon emissions that are incorporated in producing, transporting and disposing the devices, while $co_2^{op}(\mathbf{d}, \mathbf{s})$ are the operational carbon emissions. The detailed formulas are given in the appendix.
In the model, all variables are continuous and all expressions are linear. 

As future price developments are very unpredictable, we introduce an uncertainty set and consideration of regret in the model, which we will then treat in a robust way. In the following, we use the notation 
\begin{align*}
	p &= \begin{pmatrix}
		p^e \\ p^{el} \\ p^w \\p^g \\p^h
	\end{pmatrix}, 
\end{align*}
where $p^e, p^{el}, p^h, p^w, p^g \; [\euro{}/kWh]$ are the prices of buying electricity, selling electricity, buying district heat, buying wood pellets and buying gas. We define the cost regret
\begin{align*}
	costregret^{p}(\mathbf{d},\mathbf{s}) &:= costs^p(\mathbf{d},\mathbf{s})  - \min_{(\mathbf{d}^*, \mathbf{s}^*)} costs^p(\mathbf{d}^*,\mathbf{s}^*)
\end{align*}
as the difference between the minimal costs given purchasing decision $\mathbf{d}$ and price scenario $p$ and the minimal costs that are possible for this price scenario. Note that the decision variables of the minimization have to satisfy the linear constraints $A\mathbf{d}^*+B\mathbf{s}^* \leq e$ for suitably chosen parameter matrices $A,B$ and a vector $\rhs $ to represent the constraints. The detailed formulation of these constraints is given in the appendix. We obtain the model

\begin{align}
	\tag{$P(\alpha)$}\label{Palpha}
    \begin{split}
	\min_{\mathbf{d}} \max_{p} \min_{\mathbf{s}} &\begin{pmatrix}
		costregret^p(\mathbf{d},\mathbf{s}) \\ co_2(\mathbf{d},\mathbf{s}) \end{pmatrix}  \\
	\text{s.t. } &A\mathbf{d} + B\mathbf{s} \leq \rhs \\
	&A\mathbf{d}^* + B\mathbf{s}^* \leq \rhs \\
	&p \in \mathcal{P(\alpha)}.
    \end{split}
\end{align}

$\alpha$ is a parameter for the size of the price uncertainty set $\mathcal{P}(\alpha)$, for which different realizations like a box, polyhedron or ellipsoid are possible choices. 

The trade-off between $\alpha$ and solutions of \eqref{Palpha} is exactly what is addressed by inverse robustness or information gap decision theory.

\section{Solution Strategy} \label{maro}
The proposed problem formulation can be understood as a MARO problem. We follow the approach of constraint MARO-efficiency from \citep{halsermaro}, which is the application of the $\varepsilon$-constraint method in a first-scalarize-then-robustify manner to find solutions to \eqref{Palpha}. This is done by restricting the annual carbon emissions to some value $c \in \mathbb{R}$. The resulting optimization problem is

\begin{align}
    \tag{$P(\alpha, c)$} \label{Palphac}
    \begin{split}
    \min_{\mathbf{d}} \max_{p} \min_{\mathbf{s}} \;  &
	costregret^{p}(\mathbf{d},\mathbf{s})
	\\
	\text{s.t. } &  co_2(\mathbf{d},\mathbf{s}) \leq c \\
	&A\mathbf{d} + B\mathbf{s} \leq \rhs \\
	&  co_2(\mathbf{d}^*,\mathbf{s}^*) \leq c \\
	&A\mathbf{d}^* + B\mathbf{s}^* \leq \rhs \\
	&p \in \mathcal{P(\alpha)}.
    \end{split}
\end{align}

The problem contains different stages, where we call the minimization of $\mathbf{d}$ the first stage, the maximization the middle stage and the minimization over the controls $\mathbf{s}$ the second stage problem. We use this notation, where the middle stage is not counted, to be coherent with the two-stage decision structure of the problem.
To guarantee that the solution is meaningful and the algorithm from \ref{solutionalgorithm} is applicable, we have to guarantee that for every fixed $\mathbf{d}$ and all $p \in \mathcal{P}(\alpha)$ there is a feasible second stage solution $\mathbf{s}$ to \eqref{Palphac}. This assumption is called \emph{relatively complete recourse} assumption \citep{birge2011introduction}. In our case, this means that it has to be guaranteed that all loads can be fulfilled, while respecting all constraints, independent of the first stage device selection. To ensure this, we restrict the investment carbon emissions to $c$ by adding the constraint
\begin{align*}
	co_2^{inv}(\mathbf{d}) \leq c \\
\end{align*}
to the first stage problem. Further, we add dummy generators to our model that have neither investment costs nor carbon emissions, but enormous costs for generating heating and cooling. These dummy generators can always step in for providing the required energy. This guarantees relatively complete recourse for all $c > 0$.

We call solutions $\mathbf{d}$ to \eqref{Palphac} constraint MARO-efficient for \eqref{Palpha}. Moreover, for all price scenarios $p \in \mathcal{P}(\alpha)$, there is the guarantee that there is a control $\mathbf{s}$ such that for the optimal value $o$ of \eqref{Palphac}

\begin{align*}
	costregret^p(\mathbf{d}, \mathbf{s}) &\leq o \\
	co_2(\mathbf{d}, \mathbf{s}) &\leq c.
\end{align*}

To find solutions to \eqref{Palphac}, we essentially dualize the second stage problem and then apply adaptive discretization (see e.g. \citep{blankenship1976infinitely}). The detailed algorithm, including its convergence proof, is given in \ref{solutionalgorithm}. 

\section{Case Study} \label{casestudy}

\subsection{Data Basis}

For the case study, we consider a new office building in Kaiserslautern, Germany, planned next to an existing office complex. In addition to the offices, the building will also house data centers, that require significant cooling. We use the data basis from \citep{halserpn}, where the building, its heating and cooling loads, the weather conditions of the example days, the unit properties, i.e. their coefficients of performance, and the nominal price scenario are described in detail. 

\subsection{Uncertainty Set}
\begin{figure}
	\parbox{0.5\textwidth}{
		\includegraphics[width=0.5\textwidth]{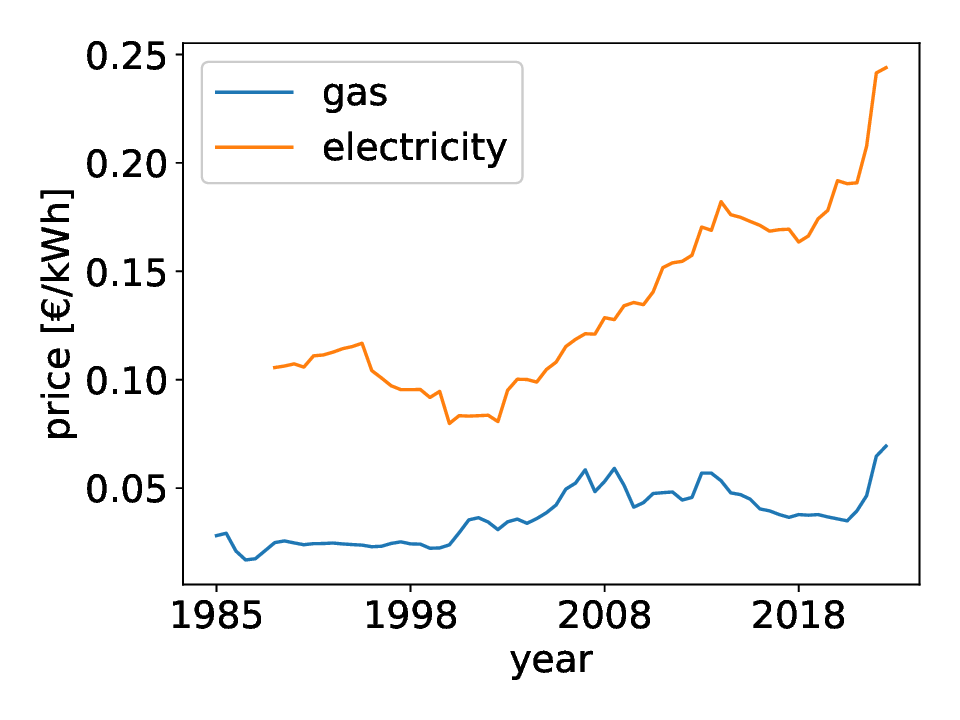}
	}
	\begin{minipage}{0.5\textwidth}
		\includegraphics[width=\textwidth]{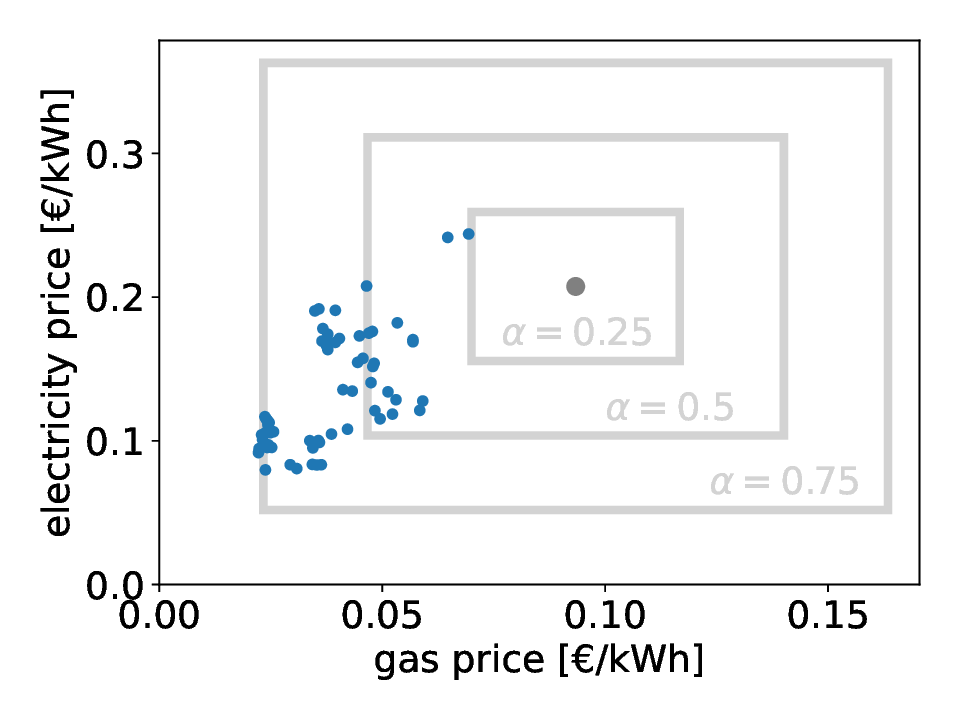}
	\end{minipage}
	\caption{Developments of gas and electricity prices for Germany from 1985 to 2022, including all taxes. Depicted are the available values for commercial prices which are closest to the prices for companies with 4 GWh of annual electricity as well as gas consumption. On the right side, the blue dots correspond to the values from the left chart, while the gray dot is a value from 2023. The gray boxes are visualizations of our choice of uncertainty set shape.}
	\label{uncertaintyset}
\end{figure}

To get a first impression of the uncertainty, we look at data from \cite{europa}, and observe the developments of the electricity and gas price from Figure~\ref{uncertaintyset} to which we added the price from 2023 (gray dot), that was assumed in \citep{halserpn}. We choose a box shape of the uncertainty set. As the values from the last year are extreme for the known observations of the uncertain parameters and future prices are unpredictable due to political and economical uncertainties, this ensures that we are likely secured against price developments of different kinds.

\begin{table}
	\centering
	\begin{tabular}{l|rrrrr}
		& $\hat{p}_{gas}$& $\hat{p}_{el_{buy}}$& $ \hat{p}_{el_{sell}}$& $\hat{p}_{dh}$ &$\hat{p}_{wp}$ \\
		\hline
		Value [\euro/kWh]& 0.0934  &  0.2074  & 0.2074 & 0.095 & 0.0817 \\
	\end{tabular}
	\caption{Nominal prices for gas, buying electricity, selling electricity, district heating and wood pellets. }\label{nominalprices}
\end{table}

We do not only consider uncertainty of gas and electricity prices but of all prices, where we assume as nominal prices the values from Table~\ref{nominalprices}. 
Then, we can define a polyhedral box uncertainty set as

\begin{align*}
	\mathcal{P}_{box}(\alpha):=& \{ (p_{gas}, p_{el_{buy}}, p_{el_{sell}},  p_{dh} , p_{wp})^T \mid \\
	& p_{gas} =\hat{p}_{gas} \cdot \beta_{gas}, p_{el_{buy}}=\hat{p}_{el_{buy}} \cdot \beta_{el_{buy}} , \\& p_{el_{sell}}= \hat{p}_{el_{sell}} \cdot \beta_{el_{sell}}, p_{dh}=\hat{p}_{dh} \cdot \beta_{dh}, \\& p_{wp}=\hat{p}_{wp} \cdot \beta_{wp}, \\ & \beta_{gas}, \beta_{el_{buy}}, \beta_{el_{sell}}, \beta_{dh}, \beta_{wp} \in [1- \alpha, 1+\alpha] \} . \\
\end{align*}

\subsection{Implementation and Computational Performance}
All computations discussed in the following were performed on an Intel i7 with 32GB RAM and operating system Windows 11. The code was written in python 3.9, using gurobipy 11.0.0.
As solution algorithms, we implemented the constraint generation (CG) and the column and constraint generation (C\&CG) algorithm, which are the two variants of Algorithm~\ref{alg:cg} that can be found in the appendix. Remark~\ref{primalcon} in the appendix suggests that using the C\&CG algorithm is superior to using the CG Algorithm.
First numerical experiments confirm this suggestion. 
Table~\ref{comptables} shows numerical results for the two algorithms for different model sizes and parameters. In the table, the respective ratios are below one for all comparisons, which means that the C\&CG algorithm has less computational time and fewer iterations for all considered scenarios. Therefore, in the following, we use the C\&CG algorithm. However, also for the C\&CG algorithm computation times increase rapidly with model size. Therefore, the problem size that is tractable is limited. Furthermore, it can be seen that not only model size but also $\alpha$ and the carbon limit influence the computational complexity. For this and all following applications of Algorithm~\ref{alg:cg}, the termination accuracy $\varepsilon$ of the Algorithm was set to 100 \euro. 

\begin{table}
	\centering
	\adjustbox{max width=0.75\textwidth}{
\begin{tabular}{l|rrrr}
	\toprule
	\multicolumn{1}{|c|}{$\alpha$} & \multicolumn{2}{c|}{0.3} & \multicolumn{2}{c|}{0.7}\\
	\cline{2-5}
	\multicolumn{1}{|c|}{Carbon Limit [t CO\textsubscript{2}]} & \multicolumn{1}{c|}{30} & \multicolumn{1}{c|}{60} & \multicolumn{1}{c|}{30} & \multicolumn{1}{c|}{60} \\
	\bottomrule
\toprule
\multicolumn{1}{c|}{$|\mathcal{U}|=1, n = 1$} & \multicolumn{4}{c}{ } \\
\midrule
\#iterations CG   &        11 &        21 &        16 &        24 \\
\#iterations C\&CG &         1 &         3 &         2 &         3 \\
\textbf{\#iterations ratio}  &         \textbf{0.09} &         \textbf{0.14} &        \textbf{ 0.12} &        \textbf{ 0.12} \\
time CG [s]         &         2.51 &         5.26 &         3.64 &         5.32 \\
time C\&CG [s]       &         0.51 &         1.30 &         0.76 &         1.19 \\
\textbf{time ratio }        &         \textbf{0.20} &         \textbf{0.25} &         \textbf{0.21} &         \textbf{0.22} \\
	\bottomrule
	\toprule
	\multicolumn{1}{c|}{$|\mathcal{U}|=1, n = 3$} & \multicolumn{4}{c}{ } \\
	\midrule
\#iterations CG   &        32 &        29 &        39 &        36 \\
\#iterations C\&CG &         1 &         3 &         2 &         3 \\
\textbf{\#iterations ratio}  &         \textbf{0.03} &         \textbf{0.10} &         \textbf{0.05} &         \textbf{0.08} \\
time CG [s]         &        18.73 &        17.61 &        44.65 &        35.48 \\
time C\&CG [s]        &         1.20 &         2.90 &         1.87 &         4.52 \\
\textbf{time ratio}         &         \textbf{0.06} &         \textbf{0.16} &         \textbf{0.04} &         \textbf{0.13} \\
	\bottomrule
\toprule
\multicolumn{1}{c|}{$|\mathcal{U}|=3, n = 1$} & \multicolumn{4}{c}{ } \\
\midrule
\#iterations CG   &        69 &        62 &        46 &        55 \\
\#iterations C\&CG &         2 &         1 &         2 &         2 \\
\textbf{\#iterations ratio}  &        \textbf{ 0.03} &         \textbf{0.02} &         \textbf{0.04} &         \textbf{0.04} \\
time CG [s]          &        63.76 &        58.23 &        45.96 &        57.38 \\
time C\&CG [s]        &         3.33 &         1.58 &         2.88 &         2.78 \\
\textbf{time ratio}         &         \textbf{0.05 }&        \textbf{ 0.03} &         \textbf{0.06} &         \textbf{0.05} \\
	\bottomrule
\toprule
\multicolumn{1}{c|}{$|\mathcal{U}|=3, n = 3$} & \multicolumn{4}{c}{ } \\
\midrule
\#iterations CG   &       141 &       164 &       102 &       137 \\
\#iterations C\&CG &         2 &         1 &         3 &         2 \\
\textbf{\#iterations ratio}  &         \textbf{0.01} &         \textbf{0.01} &         \textbf{0.03} &         \textbf{0.01} \\
time CG [s]          &       312.94 &       399.81 &       174.32 &       203.88 \\
time C\&CG [s]        &         8.82 &         3.27 &        14.89 &         5.28 \\
\textbf{time ratio}         &         \textbf{0.03} &         \textbf{0.01} &         \textbf{0.09} &         \textbf{0.03} \\
	\bottomrule
\end{tabular}
}
\caption{Comparison of CG and C\&CG algorithm for 16 different scenarios. $|\mathcal{U}|$ is the number of considered type days and $n$ is the number of time steps per day. The ratio is the ratio of the C\&CG value to the CG value.} \label{comptables}
\end{table}

\subsection{Results}
In this section, we discuss possible trade-offs between cost regret, carbon emissions and willingness to take risks, represented by $\alpha$, and their implications to device selection. 
To account for model size limitations and the chosen model formulation, we look at three different case study configurations.

First, we consider an example with three representative days in a resolution of three time steps of eight hours each and with annual carbon limits between 20 and 160 t CO\textsubscript{2} in steps of 10 t CO\textsubscript{2}. Apart from the carbon limit, we also vary the size of the uncertainty set, for which we test the $\alpha$-values of 0.25, 0.5 and 0.75.

\begin{figure}
	\centering
	\includegraphics[width=\textwidth]{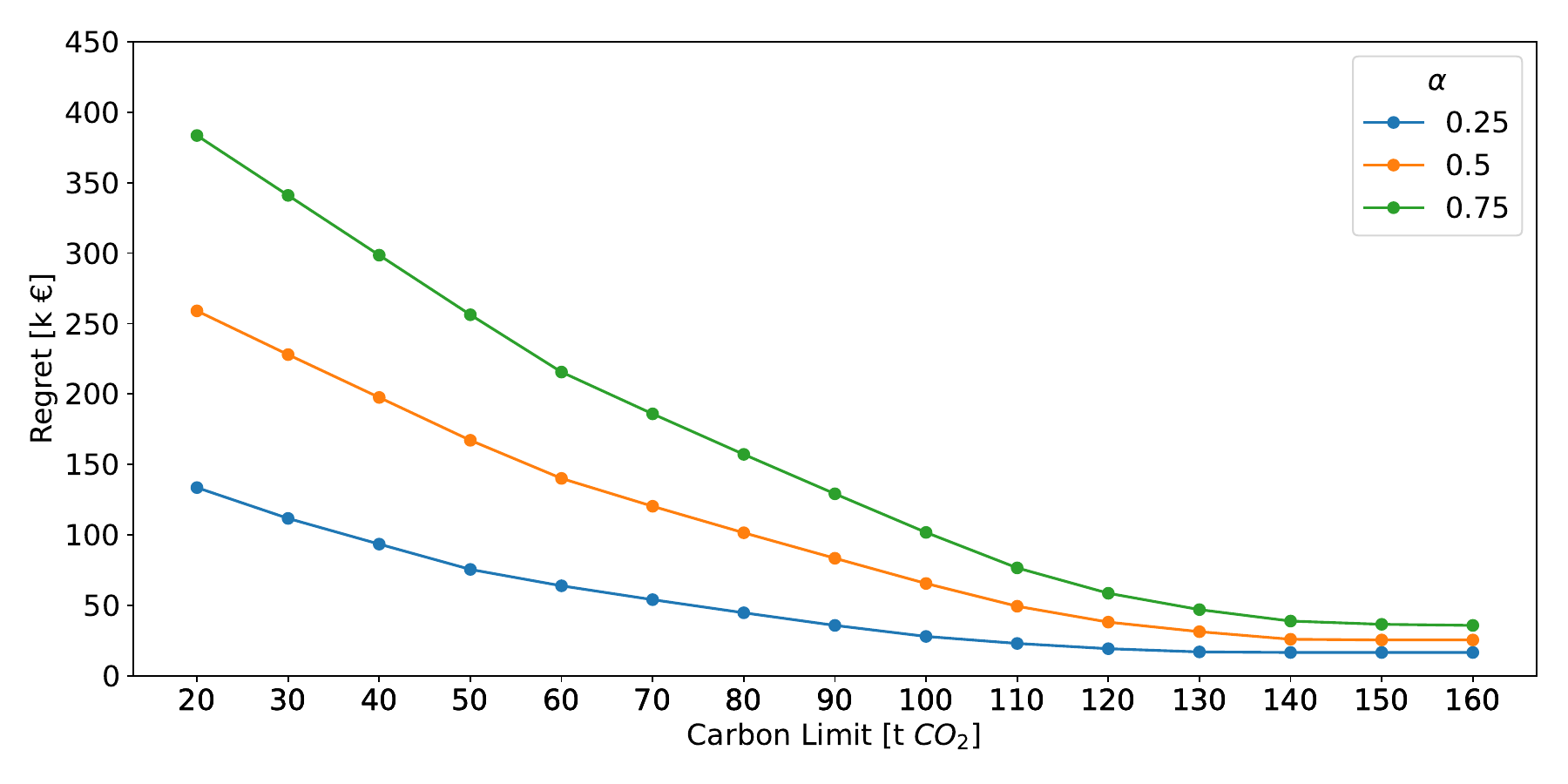}
	\includegraphics[width=\textwidth]{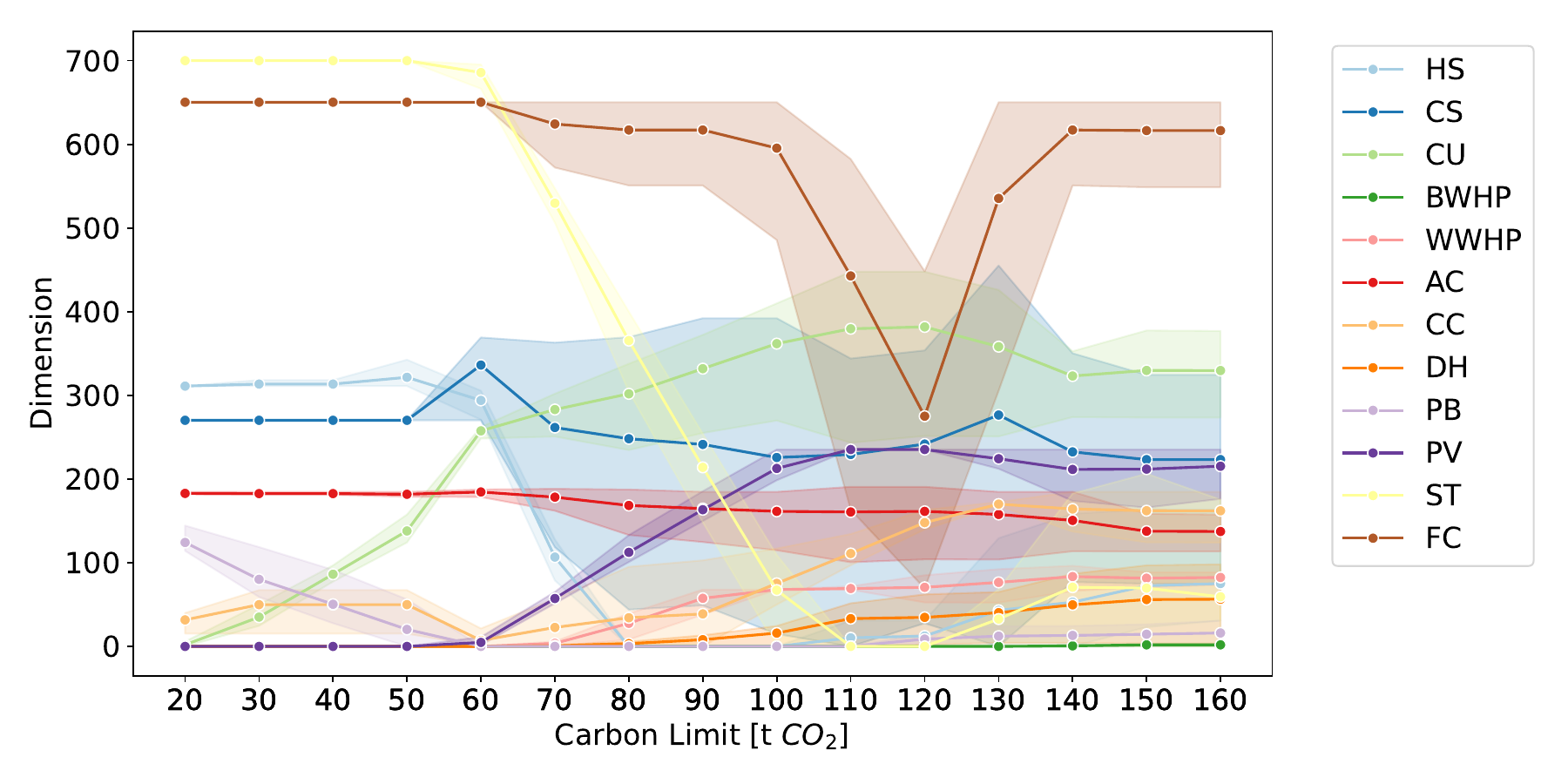}
	\caption{Overview of regret and device selection developments for different $\alpha$ and carbon limit values for three representative days consisting of three time steps, representing eight hours each. The minimal possible carbon emission value is 17170 kg, if all loads have to be fulfilled without dummy usage, as was computed by a single criteria optimization. In the lower plot, the shaded areas are the range that is caused by the different $\alpha$-values and the lines correspond to the mean of the three $\alpha$-values.}\label{results33}
\end{figure}

The visualization of the results can be seen in Figure~\ref{results33}, where in the upper plot for the different $\alpha$-values and carbon limits the cost regret of the constraint MARO-efficient solutions is shown. The connecting lines between the different carbon values are only for visual readability but not to guarantee behavior for computations with carbon limits in between. In this graphic, it is to be expected from the mathematical formulas that the regret has to be monotonically decreasing for each $\alpha$ and that higher $\alpha$-values cause higher regret, but this is the only structural knowledge that is available. The plot can be used to find a preferred trade-off, and it can be seen that for the two objectives as well as for the willingness to take risks, represented by $\alpha$, there is significant potential for trade-off finding.

In the lower plot, the influence of the set carbon limit on the selected devices is visualized. The abbreviations used in the graphic are from Table~\ref{Table}, where all possible converters with their respective dimension units are listed. Although the developments are not necessarily monotone, some tendencies can be observed. Higher carbon emission acceptance leads to smaller heat storage and solar thermal, while the size of the cogeneration unit, the water-water heat pump, district heating, compression chiller and photovoltaic is increased. Note that photovoltaic and solar thermal compete for the area on the roof. It is also important to notice that the behavior of the converter sizes does not have to be monotone over increasing carbon limits, as, for example, the value for the size of free cooling shows.

\pgfkeys{/pgf/number format/.cd,relative*={0},precision=0,relative style=fixed}
\begin{table}
	\centering
	\adjustbox{angle=0}{
		\begin{tabular}{>{\RaggedRight}l l l} 
			\hline
			converter/storage & abbreviation & dimension unit \\
			\hline
   			air-water heat pump & AWHP & kW output \\ 
   			absorption chiller & AdC & kW output \\ 
      		adsorption chiller & AC & kW output \\ 
        	brine-water heat pump & BWHP & kW output \\ 
			cogeneration unit & CU & kW output \\ 
   			cold storage & CS & kWh storage \\
   			compression chiller & CC & kW output \\ 
      		cooling dummy & C-Dummy & kW output \\
			district heating & DH & kW output \\ 
   			free cooling & FC & kW output \\ 
			gas boiler & GB &  kW output \\ 
   			heating dummy & H-Dummy & kW output \\
      		heat storage & HS & kWh storage \\
   			pellet boiler & PB & kW output \\ 
			photovoltaic & PV & $m^2$ \\ 
			reversible air-water heat pump & rev. AWHP & kW output \\ 
			reversible brine-water heat pump & rev. BWHP & kW output \\ 
			reversible water-water heat pump & rev. WWHP & kW output \\ 
			solar thermal & ST & $m^2$ \\ 
			water-water heat pump & WWHP & kW output \\ 
			\hline
		\end{tabular}
	}
	\caption{Parameters of all considered converters and storages. Note that $1m^2$ PV uses $2.975m^2$ roof area.} \label{Table}
\end{table}

Second, as three days with three time step resolution might cause a very limited case study result, we also performed computations for six days with 12 time steps of two hours each. To limit computational time, we only looked at three different carbon limits (20, 40, 60 t CO\textsubscript{2}) in this case and again at the three different $\alpha$-values of 0.25, 0.5 and 0.75. The resulting regret is compared with that of the coarser computation from the first study in Table~\ref{comparisontable}. It can be seen that despite the difference in model resolution, the regret values differ by at most 8\%. This can be seen as an indicator that the results from the coarser resolution are already quite meaningful. Therefore, it is possible to select a desired compromise in Figure~\ref{results33} and then compute the exact device configuration in a second step, with a finer resolution.

\begin{table}
	\begin{tabular}{l l c c c}
		$\alpha$ & carbon limit [t $CO_2$] & coarse regret [€] & fine regret [€] & deviation \\
		\hline
		\multirow{3}{*}{0.25}& 20 & 132801 & 133537 & 0.01 \\
		& 40 & 92615 & 93410 & 0.01 \\
		& 60 & 58487 & 63813 & 0.08 \\
		\hline
		\multirow{3}{*}{0.5}& 20 & 252957 & 258994 & 0.02\\
		& 40 & 191580 & 197533 & 0.03\\
		& 60 & 133167 & 140059 & 0.05\\
		\hline
		\multirow{3}{*}{0.75}& 20 & 374133 & 383486 & 0.02 \\
		& 40 & 288910 & 298527 & 0.03\\
		& 60 & 206311 & 215559 & 0.04\\
	\end{tabular}
	\caption{Comparison of the regret values of coarser optimization (based on three days with three time steps) and finer optimization (based on six days with twelve time steps). The deviation is the relative difference of the coarser regret compared to the finer regret.}\label{comparisontable}
\end{table}

\begin{figure}
	\centering
	\includegraphics[width=\textwidth]{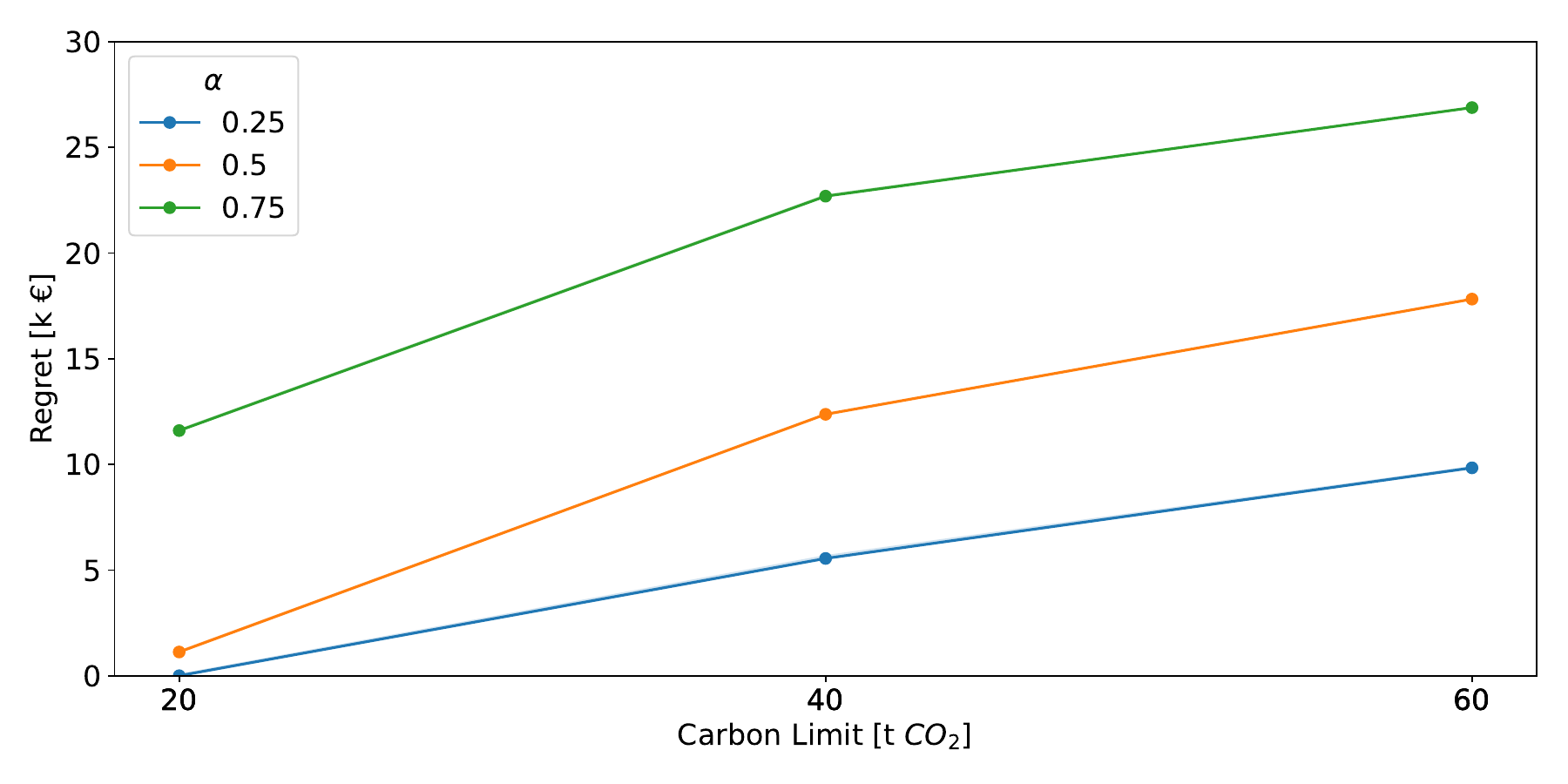}
	\caption{Overview of regret developments for different $\alpha$ and carbon limit values for three representative days consisting of three time steps, representing eight hours each. In contrast to Figure~\ref{results33}, the regret is computed in comparison with only those designs and operations that also fulfill the carbon limit.}\label{results}
\end{figure}

Third, we look at a modified version of the problem. So far, we have considered the regret when comparing with all available designs and controls. This approach would be the best choice to find trade-offs with a decision maker who considers carbon emission reduction as an optional criterion and will always ask how much money could have been saved by ignoring carbon restrictions. In the following, we will only compare to designs and controls that also respect the carbon limit. This addresses a decision maker who is interested in guaranteeing good carbon limits.
Formally, this is obtained by adding the constraints 
\begin{align*}
	co_2(\mathbf{d}^*, \mathbf{s}^*) &\leq c \\
	co_2^{inv}(\mathbf{d}^*) &\leq c
\end{align*} 
to \eqref{Palphac}.
The resulting problem structure is different from that of the so far considered constraint MARO-efficiency.
As computation times increase in this case, we again only look at three days with three time steps each and use the carbon and $\alpha$-values from the second study. The results of this study can be seen in Figure~\ref{results}. It can be observed that in this case carbon emission reduction and regret reduction are no contradicting goals anymore, but minimal carbon emissions are well in line with minimal regret. Dealing with larger uncertainty sets, by considering larger $\alpha$-values, still leads to larger regret. Considering the solutions in more detail, it can be observed that the regret-minimal configuration for this plot is for all $\alpha$-values based on solar thermal, absorption chiller, pellet boiler, cold storage, heat storage, free cooling and a small compression chiller as well as photovoltaic area. All other devices are not chosen in this case. 

\section{Conclusion} \label{conclusion}
We introduced a model and computed constraint MARO-efficient solutions for building energy supply design with the objectives of minimizing cost regret and carbon emissions under price uncertainty. Moreover, we visualized the results similarly as it is usually done for MCO, offering intuitive access to trade-off finding. For the computations, we successfully adapted an algorithm for regret with respect to load uncertainty to regret with respect to price uncertainty and proved the convergence in this case for the CG as well as for the C\&CG algorithm. 

In the case study, we made the following observations.
\begin{itemize}
	\item There is a significant potential for trade-offs between cost regret, carbon emissions and the size of the uncertainty set, which corresponds to willingness to take risks.
	\item Similar to the algorithm for load uncertainty, the problem size that is tractable is limited.
	\item The computational effort does not only depend on the model size, i.e. the number of considered days and time steps, but also on the volume of the uncertainty set and the value of the imposed carbon limit.
	\item Minimizing carbon emissions and cost regret does not lead to contradicting results if the comparison is only with devices and controls that also fulfill the carbon limit. In this case, the most carbon saving is also the most regret reducing design.
\end{itemize}

In further research, it would be beneficial to find possibilities to accelerate the algorithm to be able to deal with larger models. Moreover, it could be useful to combine the approach with more interactive decision aiding approaches like Pareto navigation, which was already applied to the problem without robustness in \citep{halserpn}. However, an additional challenge for Pareto navigation is the nonconvexity of the robust problem. Furthermore, it could be beneficial to consider additional objectives like, for example, minimization of investment costs.

\appendix

\section{Problem Formulations} \label{problemformulation}

In this section, we introduce the detailed terms of the problem formulation.

\begin{table}
	\adjustbox{max width=0.9\textwidth}{
	\begin{tabularx}{\textwidth}{p{0.23\textwidth}X}
		\toprule
		$\mathcal{M}$ & set of all possible converters \\     
  		$\mathcal{M}^d$ & set consisting of heating and cooling dummy \\
		$\mathcal{M}^{d+}$ & $\mathcal{M}$ with additional dummy generators for heating and cooling\\
  		$\mathcal{T} := \{0,...,n\}$ & set of all time steps (equidistant) of the day, where we require that n divides 24 \\
		$\mathcal{T}^+:= \mathcal{T} \cup \{n+1\}$ & set of all time steps of the day with an additional step at the end\\
		$\mathcal{U}$ & set of representative days, consisting of finite number of cluster center days, which are found with a $k$-Medoids algorithm \\
    \hline
		$\mathbf{d}_i$ & dimension of converter $i$ [-]\\       
  		$\mathbf{d}^{hs}/\mathbf{d}^{cs}$ & dimension of heat/cold storage [kWh storage potential] \\
		$\mathbf{s}_{itk}$ & dimension-scaled load of converter $i$ in time step $t$ at day $k$ [-] \\
		$\mathbf{s}_{tk}^h/\mathbf{s}_{tk}^c$ & heat/cold storage charge in time step $t$ of day $k$ [kWh] \\
    \hline
  		$c^e/c^{el}/c^h/c^w/c^g$ & carbon equivalent emissions of buying electricity/selling electricity/district heat/wood pellet/gas [kg/kWh]\\
  		$c^{hs}/c^{cs}$ & depreciated carbon emission factor of heat/cold storage [$kg$] \\
		$e_{itk}/k_{itk}/w_{itk}/$ $g_{itk}/t_{itk}$ & maximum electric/heat/wood pellet/gas/district heat consumption of converter $i$ in time step $t$ of day $k$ [kWh]\\
		$el_{itk}/h_{itk}/c_{itk}$ & maximum electric/heat/cold production of converter $i$ in time step $t$ of day $k$ [kWh]\\
  		$f := \frac{24}{n}$ & length of the time steps \\
    	$H_{tk}/C_{tk}$ & heating/cooling loads in time step $t$ of day $k$ [kWh]\\
		$p^e/p^{el}/p^h/p^w/p^g$ & price of buying electricity/selling electricity/district heat/wood pellet/gas [$\euro{}/kWh$]\\ 
        $p$ & vector of all resource prices (variable or fixed) \\
		$p_i^d/c_i^d$ & size dependent depreciation costs/carbon emission of converter $i$ [$\euro$]/[kg]\\
		$p_i^f/c_i^f$ & base depreciation costs/carbon emission of converter $i$ [$\euro$]/[kg]\\
		$p^{hs}/p^{cs}$ & depreciated price factor of heat/cold storage [$\euro$] \\
  		$w_k$ & weights of cluster day $k \in \mathcal{U}$ \\
		\bottomrule
	\end{tabularx}}
	\caption{ Notation for the LP. Indices are subscripts, while superscripts are used for nomenclature.} \label{Notation}
\end{table}

Using the notation from Table~\ref{Notation} and defining
\begin{align*}
	p_{itk}&:= e_{itk}p^e  - el_{itk}p^{el} + w_{itk}p^w + g_{itk}p^g  + t_{itk}p^h & \forall i \in \mathcal{M}^{d+} ,t \in \mathcal{T}, k \in \mathcal{U}\\
	q_{itk}&:= e_{itk}c^e - el_{itk}c^{el} + w_{itk}c^w + g_{itk}c^g  + t_{itk}c^h & \forall i \in \mathcal{M}^{d+} ,t \in \mathcal{T}, k \in \mathcal{U}\\
	p^{inv}&:= f\mathbf{d}^{hs}p^{hs} + f\mathbf{d}^{cs}p^{cs} +
	\sum_{i \in \mathcal{M}} \mathbf{d}_ip_i^d  \\
	c^{inv}&:= f\mathbf{d}^{hs}c_{hs} + f\mathbf{d}^{cs}c_{cs} +
	\sum_{i \in \mathcal{M}} \mathbf{d}_ic_i^d  \\
	\mathbf{d} &:= \begin{pmatrix}
		\mathbf{d}^{hs} \\ \mathbf{d}^{cs} \\ (\mathbf{d}_i)_{i \in \mathcal{M}^{d+}}
	\end{pmatrix} \\
	\mathbf{s} &:= \begin{pmatrix}
		(\mathbf{s}_{itk})_{i \in \mathcal{M}^{d+}, t \in \mathcal{T}, k \in \mathcal{U}} \\  (\mathbf{s}_{tk}^c)_{t \in \mathcal{T}^+, k \in \mathcal{U}}\\  
		(\mathbf{s}_{tk}^h)_{t \in \mathcal{T}^+, k \in \mathcal{U}}
	\end{pmatrix},
\end{align*}
we define the objective functions as
\begin{align*}
	costs^p(\mathbf{d},\mathbf{s})&:=f\sum_{\substack{i \in \mathcal{M}, t \in \mathcal{T}, \\ k \in \mathcal{U}}}w_k\mathbf{s}_{itk}p_{itk} + 
	p^{inv} \\
	co_2(\mathbf{d},\mathbf{s}) &:= f\sum_{\substack{i \in \mathcal{M}, t \in \mathcal{T}, \\ k \in \mathcal{U}}}w_k\mathbf{s}_{itk}q_{itk} + c^{inv}.
\end{align*}
The carbon emission split in operational and investment emissions is then given by
\begin{align*}
	co_2^{op}(\mathbf{d})&=f\sum_{\substack{i \in \mathcal{M}, t \in \mathcal{T}, \\ k \in \mathcal{U}}}w_k\mathbf{s}_{itk}q_{itk} \\ 
 co_2^{inv}(\mathbf{d}, \mathbf{s})&=c^{inv}
\end{align*}
and analogously for costs.
The detailed constraints are given as cooling and heating load fulfilling
\begin{align*}
	\mathbf{s}_{tk}^c-\mathbf{s}_{t+1k}^c + \sum_{i \in \mathcal{M}^{d+}} \mathbf{s}_{itk} c_{itk}  &= C_{tk} &&\forall t  \in \mathcal{T}, k \in \mathcal{U} \\
	\mathbf{s}_{tk}^h-\mathbf{s}_{t+1k}^h + \sum_{i \in \mathcal{M}^{d+}} \mathbf{s}_{itk} (h_{itk}-t_{itk})  &= H_{tk} &&\forall t \in \mathcal{T}, k \in \mathcal{U},\\
\end{align*}	
technical limits of the converters and storages
\begin{align*}
	\mathbf{s}_{itk} &\leq \mathbf{d}_i & \forall i \in \mathcal{M}^{d+} ,t \in \mathcal{T}, k \in \mathcal{U}\\
	\mathbf{s}_{tk}^c &\leq \mathbf{d}^{cs} &\forall t  \in \mathcal{T},  k \in \mathcal{U}\\
	\mathbf{s}_{tk}^h &\leq \mathbf{d}^{hs} &\forall t  \in \mathcal{T}, k \in \mathcal{U} \\
\end{align*}
and start and end states of the storages for each day. These are chosen such that the heat storage has to be empty at midnight and the cold storage has to be full at that time
\begin{align*}
	\mathbf{s}_{0k}^c & = 0 && \forall k \in \mathcal{U}\\
	\mathbf{s}_{\{n+1\}k}^c &= 0 && \forall k \in \mathcal{U}\\
	\mathbf{s}_{0k}^h & = \mathbf{d}^{hs} && \forall k \in \mathcal{U}\\
	\mathbf{s}_{\{n+1\}k}^h &= \mathbf{d}^{hs} && \forall k \in \mathcal{U}.\\
\end{align*}
Moreover, there are the constraints that guarantee nonnegativity of dimensions and controls and the fixation of the size of dummy generators
\begin{align*}
	p,\mathbf{s} &\geq 0  \\
	\mathbf{d} & \geq 0  \\
	\mathbf{d}_i &= M & \forall i \in \mathcal{M}^d
\end{align*}
for a sufficiently large constant $M$.
There are also device specific constraints. To avoid introducing too many indices, we give them only verbally. 
\begin{itemize}
	\item We model reversible heat pumps as two devices, where one represents the heating and the other the cooling, and both devices have to have a fixed size ratio. We ensure that only one of them has positive coefficient of performance in every time step.
	\item The roof area is an upper bound for the sum of the (weighted) dimensions of photovoltaic and solar thermal.
	\item An adsorption chiller needs to have a dimension smaller than the sum of the dimensions of the district heating and the cogeneration unit.
\end{itemize}

\section{Solution Algorithm} \label{solutionalgorithm}

\subsection{Motivation}
In the following, we motivate the solution algorithm. As this algorithm can also be applied to other problems with the same structure, we use general notation, that is independent of the problem-specific notation from before.
The minimax regret problem can be understood as a problem of the form

\begin{subequations}
	\begin{align*}
		\min_x \max_p \min_y \; &c^Tx + p^TAy - \min_{x^*, y^*}(c^Tx^* + p^TAy^*) \\
		s.t. \; &Bx + Cy \geq d \\
		&Bx^* + Cy^* \geq d \\
		&Ep \geq f \\
		&Gx \geq h \\
		& x, p, y, x^*, y^* \geq 0
	\end{align*}
\end{subequations}
with the variables $x \in \mathbb{R}^{n_x}, p \in \mathbb{R}^{n_p}, y \in \mathbb{R}^{n_y}, $ and the parameters $c \in \mathbb{R}^{n_x}, d \in \mathbb{R}^{n_d}, f \in \mathbb{R}^{n_f}, h \in \mathbb{R}^{n_h}, A \in \mathbb{R}^{n_p \times n_y}, B \in \mathbb{R}^{n_d \times n_x}, C \in \mathbb{R}^{n_d \times n_y}, E \in \mathbb{R}^{n_f \times n_p}, G \in \mathbb{R}^{n_h \times n_x}$ for dimensions $n_x, n_p, n_y, n_d, n_f, n_h \in \mathbb{N}$.

This problem can be reformulated as

\begin{subequations} \label{primal}
	\begin{align}
		\min_x \max_{p, x^*, y^*} \min_y \; &c^Tx + p^TAy - c^Tx^* - p^TAy^* \\
		s.t. \; &Bx + Cy \geq d \\
		&Bx^* + Cy^* \geq d \\
		&Ep \geq f \\
		&Gx \geq h \\
		& x, p, y, x^*, y^* \geq 0.
	\end{align}
\end{subequations}

By dualizing the second stage problem, we obtain with dual variable vector $\pi \in \mathbb{R}^{n_d}$

\begin{subequations} \label{dualized}
	\begin{align}
		\min_x \max_{p, x^*, y^*, \mathbf{\pi}} \; &c^Tx + (d-Bx)^T \mathbf{\pi} - c^Tx^* - p^TAy^* \\
		s.t. \; &C^T\mathbf{\pi} \leq A^Tp \\
		&Bx^* + Cy^* \geq d \\
		&Ep \geq f \\
		&Gx \geq h \\
		& x, p, \mathbf{\pi}, x^*, y^* \geq 0.
	\end{align}
\end{subequations}

We can observe that the feasible set of the second stage problem is independent of $x$ and has bilinear terms. We will exploit this structure in the following solution algorithms.

\subsection{Algorithm Formulation}
The algorithms we propose are based on the algorithms described by \cite{zeng2013solving}, with the difference that in our dualized problem $p$ appears not only in the objective function but also in the constraints. Therefore, an argument is necessary to observe that the optimal value in this case is still in a vertex of a certain polyhedron. To state the algorithms, we first have to split the problem in main- and subproblem. There are two possible algorithms, a constraint generation (CG) and a column and constraint generation (C\&CG) one, with slightly different main problem, which we look at in the following.

Let us start with the subproblem. In the following, we assume relatively complete recourse. This means that for every first and middle stage decision, there is a feasible second stage decision. We obtain the subproblem

\begin{align*}
	SP(x): \; \; sp(x) = \max_{p, x^*, y^*, \mathbf{\pi}} \; &(d-Bx)^T \mathbf{\pi} -c^T x^* -p^TAy^* \\
	s.t. \; &C^T \mathbf{\pi} \leq A^Tp \\
	&  Bx^* + Cy^* \geq d \\
	&Ep \geq f \\
	& p, \mathbf{\pi}, x^*, y^* \geq 0,
\end{align*}
which is a bilinear problem (BLP). We will observe later that the solution to this problem can be assumed to be in a vertex of the polyhedron that defines the feasible region of the problem. For the rest of the article, we assume that whenever we talk about a solution to the problem, it is chosen in this way. Similar BLPs can be linearized, for example in \citep{yokoyama2014revised}, but as the constraints involving $p$ are not independent of other variables but also contain $\mathbf{\pi}$, we cannot apply that to the structure of our problem and therefore have to keep the bilinear form.  

For the main problem, we consider with additional variable $\eta \in \mathbb{R}$ the epigraph reformulation of \eqref{dualized}
\begin{align*}
	\min_{x, \eta}  \; &c^Tx + \eta \\
	s.t. \; & \eta \geq (d-Bx)^T \mathbf{\pi} - c^Tx^* - p^TAy^* & \forall (p, x^*, y^*, \mathbf{\pi}) \in S  \\
	&Gx \geq h \\
	& x \geq 0,
\end{align*}
where 
\begin{align*}
	S:= \left\{(p, x^*, y^*, \mathbf{\pi}) \in \mathbb{R}_{\geq 0} ^{n_p+n_x+n_y+n_d} \mid 	C^T\mathbf{\pi} \leq A^Tp, \; Bx^* + Cy^* \geq d, \; Ep \geq f \right\}.
\end{align*}
Problems like this with an infinite number of constraints are call semi-infinite.
By relaxing the semi-infinite constraint to a still to be concretized, finite subset 
\begin{align*}
	S_k := \left\{(p_l, x_l^*, y_l^*, \mathbf{\pi}_l) \mid (p_l, x_l^*, y_l^*, \mathbf{\pi}_l) \in S, \; l \in \{0,...,k\}  \right\},
\end{align*}
we obtain the main problem
\begin{align*}
	MP^{CG}(k): \; \; \min_{x, \eta} \; &c^Tx + \eta \\
	s.t. \; & \eta \geq (d-Bx)^T \mathbf{\pi}_l -c^Tx_l^*-p_l^TAy_l^* & \forall l \leq k \\
	&Gx \geq h \\
	& x \geq 0, \eta \in \mathbb{R}.
\end{align*}
The algorithm that is based on this problem formulation is the CG algorithm. It is also possible to consider an epigraph reformulation of \eqref{primal}, which leads to the C\&CG algorithm with main problem
\begin{align*}
	MP^{C\&CG}(k): \; \; \min_{x, \eta, y_1,..., y_k} \; &c^Tx + \eta \\
	s.t. \; & \eta \geq p_l^TAy_l -c^Tx_l^*-p_l^TAy_l^* & \forall l \leq k \\
	&Bx + Cy_l \geq d & \forall l \leq k \\
	&Gx \geq h \\
	& x \geq 0, \eta \in \mathbb{R}.
\end{align*}

Both problems are LPs, so without loss of generality we assume in the following that solutions are obtained at vertices of the polyhedron defining the feasible region.
Now that we have stated all partial problems, we can state Algorithm~\ref{alg:cg}. Depending on whether we use $MP^{CG}$ or $MP^{C\&CG}$ therein as $MP$, we are in the case of the CG or the C\&CG algorithm.

\begin{algorithm}[H]
	\caption{(Column and) Constraint Generation}\label{alg:cg}
	\begin{algorithmic}
				\Require $\varepsilon \geq 0, lb=-\infty, ub=+\infty, k=0$
				\State $ x_0, \eta_0 \gets MP(0)$
				\State $ p_1, x_1^*, y_1^*, \mathbf{\pi}_1 \gets SP(x_0)$
				\State $lb \gets \max\{lb, c^Tx_0 + \eta_0\}$
				\State $ub \gets \min\{ub, c^Tx_0 + sp(x_0)\}$
				\State $k=1$
				\While{$ub-lb > \varepsilon$}
				\State $x_{k}, \eta_k \gets MP(k)$
				\State $p_{k+1}, x_{k+1}^*, y_{k+1}^*, \mathbf{\pi}_{k+1} \gets SP(x_k)$
				\State $lb \gets \max\{lb, c^Tx_k + \eta_k\}$
				\State $ub \gets \min\{ub, c^Tx_k + sp(x_k)\}$
				\State $k \gets k+1$
				\EndWhile \\
				\Return $x_k$
	\end{algorithmic}
\end{algorithm}

By iteratively solving main- and subproblem, we compute lower bounds ($lb$) and upper bounds ($ub$) for the problem. This is done until a predefined solution quality $\varepsilon \geq 0 $ is reached. As the solution to the bilinear subproblem can be very time-consuming, it can be beneficial in practice not to solve it to full optimality but to use upper bounds instead of solutions in the algorithm.

\subsection{Proof of Convergence}
\begin{thm} \label{dualconv}
	The CG algorithm terminates in a finite number of steps. The difference between the returned solution and the global optimal solution is at most $\varepsilon$.
\end{thm}

\begin{pf}
	The proof is similar to the proof in A 2.2. by \cite{yokoyama2014revised}. It follows from the definition of the main- and the subproblem that the respective solutions are lower and upper bounds for the problem. With the definitions
	\begin{align*}
		h_1=\begin{pmatrix}
			\mathbf{\pi} \\ p
		\end{pmatrix}, h_2= \begin{pmatrix}
			x^* \\ y^*
		\end{pmatrix}
	\end{align*}
	the subproblem is a bilinear problem
	\begin{align*}
		SP(x): \; \; \max_{h_1, h_2} \; & \begin{pmatrix}
			d-Bx & 0 
		\end{pmatrix} h_1 - \begin{pmatrix}
			c & 0 
		\end{pmatrix}h_2 - h_1^T \begin{pmatrix}
			0 & 0 \\ 0 & A
		\end{pmatrix} h_2 \\
		s.t. \; & \begin{pmatrix}
			-C^T& A^T\\
			0 & E
		\end{pmatrix}h_1 \geq \begin{pmatrix}
			0 \\ f
		\end{pmatrix} \\
		&  \begin{pmatrix}
			B & C
		\end{pmatrix}h_2 \geq d \\
		& h_1, h_2 \geq 0 
	\end{align*}
	to which Theorem 2.1 from \cite{konno1976cutting} can be applied, and therefore, the solution has to be obtained in a vertex of the polyhedron defining the feasible region, which is independent of $x$. When generating constraints in the main problem corresponding to all vertices of the polyhedron that defines the feasible region of the subproblem, the solution is obtained and the algorithm terminates as main- and subproblem solution coincide.
	
	It remains to show that the algorithm cannot run into an infinite loop where the same constraints are added to the mainproblem repeatedly, while others are never considered. In this case, there would be an iteration, where all considered constraints are added and all further iterations of the loop only add the same constraints. From this iteration on, the feasible region of the main problem will not change anymore. As we assume mainproblem solutions to be in vertices of the feasible region, there will be $x_{k-n}$ and $x_k$ such that $x_{k-n}=x_k$ for $k,n \in \mathbb{N}, k \geq n$ in a finite number of steps. But then, $x_k$ has to be the solution to the problem and the algorithm terminates.
\end{pf}

\begin{rem} \label{primalcon}
	The C\&CG algorithm converges in a finite number of steps, where the maximum possible number of steps is at most as large as the maximum possible number of steps for the dual algorithm. The convergence proof is analogous to that of Theorem~\ref{dualconv}, with the only difference that the constraints only correspond to $p$-values of the vertices of the polyhedron that defines the feasible region. As multiple vertices might share the same $p$-value, this number can be significantly smaller than the number of constraints that are required in the dual algorithm.
\end{rem}

\begin{rem}
	Both variants of Algorithm~\ref{alg:cg} also work with discrete or mixed integer first stage decision.
\end{rem}




\section*{CRediT authorship contribution statement}
E. Halser: Conceptualization, Data curation, Formal analysis, Investigation, Methodology, Visualization, Writing - original draft.
E. Finhold: Supervision, Writing - review \& editing.
N. Leithäuser: Supervision, Writing - review \& editing.
T. Seidel: Supervision, Writing - review \& editing.
K.-H. Küfer: Supervision.

\section*{Declaration of Generative AI and AI-assisted technologies in the writing process}
During the preparation of this work, the authors used FhGenie and DeepL Write in order to improve the language. After using these services, the authors reviewed and edited the content as needed and take full responsibility for the content of the publication.

\section*{Data Availability}
Data will be made available on request.

\section*{Funding Sources}
This research did not receive any specific grant from funding agencies in the public, commercial, or not-for-profit sectors.

\section*{Acknowledgements}
We thank our colleagues Jean Bertrand Gauthier for helpful discussions concerning column and constraint generation and Kerstin Schneider for formulation help while writing the article.

\bibliographystyle{elsarticle-harv} 
\bibliography{all_sources}


%
%
%
\end{document}